\documentclass{amsart}

\newtheorem{theorem}{Theorem}[section]
\newtheorem{lemma}[theorem]{Lemma}
\newtheorem{corollary}[theorem]{Corollary}
\newtheorem{proposition}[theorem]{Proposition}

\theoremstyle{definition}

\theoremstyle{remark}

\numberwithin{equation}{section}

\begin{document}

\title[Kloosterman Sums with Trace One Arguments]{    $\begin{array}{c}
  \text{An Infinite Family of Recursive Formulas Generating Power Moments}\\
 \text{of Kloosterman Sums with Trace One Arguments: $O(2n+1, 2^r)$ Case}
\end{array} $            }

\author{dae san kim}
\address{Department of Mathematics, Sogang University, Seoul 121-742, Korea}
\curraddr{Department of Mathematics, Sogang University, Seoul
121-742, Korea}
\email{dskim@sogong.ac.kr}
\thanks{}

%    General info
\subjclass[]{}

\date{}

\dedicatory{ }


\begin{abstract}
In this paper, we construct an infinite family of binary linear codes
associated with double cosets with respect to certain maximal parabolic
subgroup of the orthogonal group $O(2n+1, q)$.  Here $q$ is a power of two.
Then we obtain an infinite family of recursive fomulas generating the odd power
moments of Kloosterman sums with trace one arguments in terms of the frequencies
of weights in the codes associated with those double cosets in $O(2n+1,q)$
and in the codes associated with similar double cosets in the symplectic
group $Sp(2n,q)$.  This is done via Pless power moment identity and by
utilizing the explicit expressions of
exponential sums over those double cosets related to the evaluations of ``Gauss
sums" for the orthogonal group $O(2n+1,q)$.\\

Index terms-Kloosterman sum, trace one, orthogonal group, symplectic group,
double cosets, maximal parabolic subgroup, Pless power moment identity,
weight distribution.\\

\end{abstract}

\maketitle

%%%%%%%%%%%%%%%%%%%%%%%%%%%%%%%%%%%%%%%%%%%%%%%%%%%%%%%%%%%%%%%%%%%%%%%%
%%%%%%%%%%%%%%%%%%%%%%%%%%%%%%%%%%%%%%%%%%%%%%%%%%%%%%%%%%%%%%%%%%%%%%%%
\section{Introduction}
%%%%%%%%%%%%%%%%%%%%%%%%%%%%%%%%%%%%%%%%%%%%%%%%%%%%%%%%%%%%%%%%%%%%%%%%
%%%%%%%%%%%%%%%%%%%%%%%%%%%%%%%%%%%%%%%%%%%%%%%%%%%%%%%%%%%%%%%%%%%%%%%%

 Let $\psi$ be a nontrivial additive character of the
finite field $\mathbb{F}_{q}$ with $q=p^r$ elements ($p$ a prime).
Then the Kloosterman sum $K(\psi;a)$ (\cite{K1}) is defined as

\begin{align*}
K(\psi;a)=\sum_{\alpha\in
\mathbb{F}_{q}^{*}}\psi(\alpha_{}^{}+a\alpha_{}^{-1})~
(a\in\mathbb{F}_{q}^{*}).
\end{align*}
For this, we have the Weil bound
%(1)%%%%%%%%%%%%%%%%%%%%%%%%%%%%%%%%%%%%%%%%%%%%%%%%%%%%%%%%%%%%%%
\begin{equation}\label{a}
|K(\psi;a)|\leq 2\sqrt{q}.
\end{equation}
%%%%%%%%%%%%%%%%%%%%%%%%%%%%%%%%%%%%%%%%%%%%%%%%%%%%%%%%%%%%%%%%%%
The Kloosterman sum was introduced in 1926(\cite{DJ1}) to give an
estimate for the Fourier coefficients of modular forms.

For each nonnegative integer $h$, by $MK(\psi)^{h}$ we will denote
the $h$-th moment of the Kloosterman sum $K(\psi;a)$. Namely, it
is given by
\begin{align*}
MK(\psi)^{h}=\sum_{a\in\mathbb{F}_{q}^{*}}K(\psi;a)^{h}.
\end{align*}
If $\psi=\lambda$ is the canonical additive character of
$\mathbb{F}_{q}$, then $MK(\lambda)^{h}$ will be simply denoted by
$MK^{h}$.

Explicit computations on power moments of Kloosterman sums were
begun with the paper \cite{PV1} of Sali\'{e} in 1931, where he
showed, for any odd prime $q$,

\begin{align*}
MK^h=q^2 M_{h-1} -(q-1)^{h-1}+2(-1)^{h-1}~(h\geq1).
\end{align*}
Here $M_0=0$, and, for $h\in\mathbb{Z}_{>0}$,
\begin{align*}
M_h=|\{(\alpha_1,\cdots,\alpha_h)\in(\mathbb{F}_q^*)^h|\sum_{j=1}^h
\alpha_j=1=\sum_{j=1}^h \alpha_j^{-1}\}|.
\end{align*}
For $q=p$ odd prime, Sali\'{e} obtained $MK^1$, $MK^2$, $MK^3$,
$MK^4$ in \cite{PV1} by determining $M_1$, $M_2$, $M_3$. On the other
hand, $MK^5$ can be expressed in terms of the $p$-th eigenvalue for
a weight 3 newform on $\Gamma_0(15)$(cf. \cite{LN1}, \cite{M1}).
$MK^6$ can be expressed in terms of the $p$-th eigenvalue for a weight 4
newform on $\Gamma_0(6)$(cf. \cite{HS1}). Also, based on numerical
evidence, in \cite{GS1} Evans was led to propose a conjecture which
expresses $MK^7$ in terms of Hecke eigenvalues for a weight 3
newform on $\Gamma_0(525)$ with quartic nebentypus of conductor 105.

From now on, let us assume that $q=2^r$.  Carlitz\cite{E1} evaluated
$MK^h$ for $h\leq4$.  Recently, Moisio was able to find explicit
expressions of $MK^h$, for $h\leq10$~(cf.\cite{MS1}). This was done,
via Pless power moment identity, by connecting moments of
Kloosterman sums and the frequencies of weights in the binary Zetterberg
code of length $q+1$, which were known by the work of Schoof and Vlugt
in \cite{S1}.

In order to describe our results, we introduce two incomplete power moments of Kloosterman sums,
namely, the one with the sum over all $a$ in $\mathbb{F}_{q}^{*}$ with
$\emph{tr}~ a$=0
and the other with the sum over all $a$ in $\mathbb{F}_{q}^{*}$
with $\emph{tr}~ a=1$.
For every nonnegative integer $h$, and $\psi$ as before, we define

%(2)%%%%%%%%%%%%%%%%%%%%%%%%%%%%%%%%%%%%%%%%%%%%%%%%%%%%%%%%%%%%%%%%
\begin{equation}\label{b}
T_{0}K(\psi)^{h}=\sum_{a\in\mathbb{F}_{q}^{*},~tra=0}K(\psi;a)^h,~
 T_{1}K(\psi)^{h}=\sum_{a\in\mathbb{F}_{q}^{*},~tra=1}K(\psi;a)^h,
\end{equation}
%%%%%%%%%%%%%%%%%%%%%%%%%%%%%%%%%%%%%%%%%%%%%%%%%%%%%%%%%%%%%%%%%%%%
which will be respectively  called the $h$-th moment of Kloosterman sums with
$\lq\lq$trace zero arguments" and those with $\lq\lq$trace one arguments".
Then, clearly we have
%(3)%%%%%%%%%%%%%%%%%%%%%%%%%%%%%%%%%%%%%%%%%%%%%%%%%%%%%%%%%%%%%%%%
\begin{equation}\label{c}
MK(\psi)^{h}=T_{0}K(\psi)^{h}+T_{1}K(\psi)^{h}.
\end{equation}
%%%%%%%%%%%%%%%%%%%%%%%%%%%%%%%%%%%%%%%%%%%%%%%%%%%%%%%%%%%%%%%%%%%%

If $\psi=\lambda$ is the canonical additive character of
$\mathbb{F}_q$, then $T_{0}K(\lambda)^{h}$ and $T_{1}K(\lambda)^{h}$ will
be respectively denoted by $T_{0}K^{h}$ and $T_{1}K^{h}$,
for brevity.

In \cite{D4}, we obtained a recursive formula generating the odd power
moments of Kloosterman sums with trace one arguments. This was expressed
in terms of the frequencies of weights in the binary linear codes
$C(O(3,q))$ and $C(Sp(2,q))$, respectively associated with the orthogonal group
$O(3,q)$ and the symplectic group $Sp(2,q)$.

In this paper, we will show the main Theorem \ref{A} giving an infinite family of
recursive formulas generating the odd power moments of Kloosterman sums
with trace one arguments. To do that, we construct binary linear codes
$C(DC(n,q))$, associated with the double cosets $DC(n,q)$=$P\sigma_{n-1}P$, for
the maximal parabolic subgroup $P$=$P(2n+1,q)$ of the orthogonal group $O(2n+1,q)$,
and express those power moments in terms of the frequencies of weights in the codes
$C(DC(n,q))$ and $C(\widehat{DC}(n,q))$.  Here $C(\widehat{DC}(n,q))$
is a binary linear code constructed similarly from certain double cosets
$\widehat{DC}(n,q)$ in the sympletic group $Sp(2n,q)$.
Then, thanks to our previous results on the explicit expressions of exponential sums
over those double cosets related to the evaluations of ``Gauss sums" for the orthogonal
group $O(2n+1,q)$~\cite{D5}, we can express the weight of each codeword in the dual of
the codes $C(DC(n,q))$ in terms of Kloosterman sums. Then our formulas will follow
immediately from the Pless power moment identity.

Theorem \ref{A} in the following(cf.~(\ref{f})-(\ref{h})) is the main result of this paper.
Henceforth, we agree that the binomial coefficient $\begin{pmatrix}
                                                       b \\
                                                       a
                                                     \end{pmatrix}=0$
if $a>b$ or $a<0$. To simplify notations, we introduce the following ones which will
be used throughout this paper at various places.

%(4)%%%%%%%%%%%%%%%%%%%%%%%%%%%%%%%%%%%%%%%%%%%%%%%%%%%%%%%%%%%%%%%%
\begin{equation}\label{d}
A(n,q)=q^{\frac{1}{4}(5n^2-1)}\begin{bmatrix}
                        n \\
                        1 \\
                      \end{bmatrix}_{q}\Pi_{j=1}^{(n-1)/2}(q^{2j-1}-1),
\end{equation}
%%%%%%%%%%%%%%%%%%%%%%%%%%%%%%%%%%%%%%%%%%%%%%%%%%%%%%%%%%%%%%%%%%%%

%(5)%%%%%%%%%%%%%%%%%%%%%%%%%%%%%%%%%%%%%%%%%%%%%%%%%%%%%%%%%%%%%%%%
\begin{equation}\label{e}
B(n,q)=q^{\frac{1}{4}(n-1)^2}(q^{n}-1)\Pi_{j=1}^{(n-1)/2}(q^{2j}-1).
\end{equation}\\
%%%%%%%%%%%%%%%%%%%%%%%%%%%%%%%%%%%%%%%%%%%%%%%%%%%%%%%%%%%%%%%%%%%%

%Theorem 1%%%%%%%%%%%%%%%%%%%%%%%%%%%%%%%%%%%%%%%%%%%%%%%%%%%%%%%%%%
\begin{theorem}\label{A}
Let $q=2^{r}$. Assume that $n$ is any odd integer$\geq3$,
with all q, or n=1, with $q\geq8$.  Then, in the notations of (\ref{d})
and (\ref{e}), we have the following. For h=1,3,5,$\cdots$,

%(6)%%%%%%%%%%%%%%%%%%%%%%%%%%%%%%%%%%%%%%%%%%%%%%%%%%%%%%%%%%%%%%%%
\begin{equation}\label{f}
\begin{split}
&T_{1}K^{h}=-\sum_{0 \leq l\leq h-2, ~l~ odd}\begin{pmatrix}
                                         h \\
                                         l \\
                                       \end{pmatrix}B(n,q)^{h-l}T_{1}K^{l}\\
&+qA(n,q)^{-h}\sum_{j=0}^{min\{N(n,q),h\}}(-1)^{j}D_{j}(n,q)\sum_{t=j}^{h}t!
S(h,t)2^{h-t-1}\begin{pmatrix}
                 N(n,q)-j \\
                 N(n,q)-t \\
               \end{pmatrix},
\end{split}
\end{equation}
%%%%%%%%%%%%%%%%%%%%%%%%%%%%%%%%%%%%%%%%%%%%%%%%%%%%%%%%%%%%%%%%%%%%
where $N(n,q)=|DC(n,q)|=A(n,q)B(n,q), ~D_{j}(n,q)=C_{j}(n,q)-\widehat{C}_{j}(n,q)$,
with $\{C_{j}(n,q)\}_{j=0}^{N(n,q)}$, $\{\widehat{C}_{j}(n,q)\}_{j=0}^{N(n,q)}$
respectively the weight distributions of the binary linear codes
$C(DC(n,q))$ and $C(\widehat{DC}(n,q))$ given by: for $j=0,\cdots,N(n,q)$,

%(7)%%%%%%%%%%%%%%%%%%%%%%%%%%%%%%%%%%%%%%%%%%%%%%%%%%%%%%%%%%%%%%
\begin{equation}\label{g}
\begin{split}
C_{j}(n,q)=&\sum\begin{pmatrix}
                  q^{-1}A(n,q)(B(n,q)+1) \\
                  \nu_{1} \\
                \end{pmatrix}\\
&\times\prod_{tr(\beta-1)^{-1}=0}\begin{pmatrix}
                  q^{-1}A(n,q)(B(n,q)+q+1) \\
                  \nu_{\beta} \\
                \end{pmatrix}\\
                &\qquad\times\prod_{tr(\beta-1)^{-1}=1}\begin{pmatrix}
                  q^{-1}A(n,q)(B(n,q)-q+1) \\
                  \nu_{\beta} \\
                \end{pmatrix},
\end{split}
\end{equation}
%%%%%%%%%%%%%%%%%%%%%%%%%%%%%%%%%%%%%%%%%%%%%%%%%%%%%%%%%%%%%%%%%%
\\

%(8)%%%%%%%%%%%%%%%%%%%%%%%%%%%%%%%%%%%%%%%%%%%%%%%%%%%%%%%%%%%%%%
\begin{equation}\label{h}
\begin{split}
\widehat{C}_{j}(n,q)=&\sum\begin{pmatrix}
                  q^{-1}A(n,q)(B(n,q)+1) \\
                  \nu_{0} \\
                \end{pmatrix}\\
&\times\prod_{tr(\beta^{-1})=0}\begin{pmatrix}
                  q^{-1}A(n,q)(B(n,q)+q+1) \\
                  \nu_{\beta} \\
                \end{pmatrix}\\
                &\qquad\times\prod_{tr(\beta^{-1})=1}\begin{pmatrix}
                  q^{-1}A(n,q)(B(n,q)-q+1) \\
                  \nu_{\beta} \\
                \end{pmatrix}.
\end{split}
\end{equation}
%%%%%%%%%%%%%%%%%%%%%%%%%%%%%%%%%%%%%%%%%%%%%%%%%%%%%%%%%%%%%%%%%%%%%%%%%
\end{theorem}
%%%%%%%%%%%%%%%%%%%%%%%%%%%%%%%%%%%%%%%%%%%%%%%%%%%%%%%%%%%%%%%%%%%%%%%%%

Here the first sum in (\ref{f}) is 0 if $h=1$ and the unspecified sums in
(\ref{g}) and (\ref{h}) are over all the sets of nonnegative integers $\{{\nu_{\beta}}\}_{\beta\in
\mathbb{F}_{q}}$ satisfying $\sum_{\beta\in\mathbb{F}_{q}}\nu_{\beta}=j$ and
$\sum_{\beta\in\mathbb{F}_{q}}\nu_{\beta}\beta=0$.  In addition, $S(h,t)$
is the Stirling number of the second kind defined by

%9%%%%%%%%%%%%%%%%%%%%%%%%%%%%%%%%%%%%%%%%%%%%%%%%%%%%%%%%%%%%%%%%
\begin{equation}\label{i}
S(h,t)=\frac{1}{t!}\sum_{j=0}^{t}(-1)^{t-j}\begin{pmatrix}
                                             t \\
                                             j \\
                                           \end{pmatrix}j^{h}.
\end{equation}\\\\
%%%%%%%%%%%%%%%%%%%%%%%%%%%%%%%%%%%%%%%%%%%%%%%%%%%%%%%%%%%%%%%%%%

%%%%%%%%%%%%%%%%%%%%%%%%%%%%%%%%%%%%%%%%%%%%%%%%%%%%%%%%%%%%%%%%%%%%%%%%
%%%%%%%%%%%%%%%%%%%%%%%%%%%%%%%%%%%%%%%%%%%%%%%%%%%%%%%%%%%%%%%%%%%%%%%%
\section{$O(2n+1,q)$}
%%%%%%%%%%%%%%%%%%%%%%%%%%%%%%%%%%%%%%%%%%%%%%%%%%%%%%%%%%%%%%%%%%%%%%%%
%%%%%%%%%%%%%%%%%%%%%%%%%%%%%%%%%%%%%%%%%%%%%%%%%%%%%%%%%%%%%%%%%%%%%%%%
For more details about this section, one is referred to the paper
\cite{D5}. Throughout this paper, the following notations
will be used:
\begin{align*}
\begin{split}
q&=2^r~(r\in\mathbb{Z}_{>0}),\\
\mathbb{F}_q&=~the~finite~field~ with~ q~ elements,\\
TrA&=~the~ trace~ of~ A~ for~ a~ square~ matrix~ A,\\
^t B&=~the~transpose~ of~ B~for~any~matrix~B.
\end{split}
\end{align*}

Let $\theta$ be the nondegenerate quadratic form on the vector space
$\mathbb{F}_{q}^{(2n+1)\times 1}$ of all $(2n+1)\times 1$ column vectors over
$\mathbb{F}_{q}$, given by
\begin{align*}
\theta(\sum_{i=1}^{2n+1}x_{i}e^{i})=\sum_{i=1}^{n}x_{i}x_{n+i}+x_{2n+1}^{2},
\end{align*}
where $\{e^{1}=~^t [10\cdots0], e^{2}~=~^t [010\cdots0],\cdots,e^{2n+1}=~^t[0\cdots01]\}$
is the standard basis of $\mathbb{F}_{q}^{(2n+1)\times1}$.

The group $O(2n+1,q)$ of all isometries of $(\mathbb{F}_{q}^{(2n+1)\times1},\theta)$
consists of the matrices

\begin{align*}
\begin{bmatrix}
  A & B & 0 \\
  C & D & 0 \\
  g & h & 1 \\
\end{bmatrix}(A,B,C,D~n\times n,g,h~1\times n)
\end{align*}
in $GL(2n+1,q)$ satisfying the relations:
\begin{align*}
\begin{split}
&{}^t AC+{{}^tgg}~is~alternating\\
&{{}^t BD}+{}^t hh~is~alternating\\
&{}^t AD+{}^t CB=1_{n}.
\end{split}
\end{align*}
Here an $n\times n$ matrix $(a_{ij})$ is called alternating if
\begin{align*}
\left\{
  \begin{array}{ll}
    a_{ii}~=~0, & \hbox{$~for~1\leq i\leq n,$} \\
    a_{ij}~=~-a_{ji}~=~a_{ji}, & \hbox{$~for~ 1\leq i < j \leq n.$}
  \end{array}
\right.
\end{align*}

Also, one observes, for example, that $^t AC+{^tg}g$ is alternating
if and only if $^t AC={^t C}A$ and $g=\sqrt{diag(^t AC)}$, where
$\sqrt{diag(^t AC)}$ indicates the $1\times n$ matrix $[\alpha_{1},\cdots,\alpha_{n}]$
if the diagonal entries of $^t AC$ are given by
\begin{align*}
({^t AC})_{11}~=~\alpha_{1}^{2},\cdots,({^t AC})_{nn}~=~\alpha_{n}^{2},~for~\alpha_{i}\in
\mathbb{F}_{q}.
\end{align*}

As is well known, there is an isomorphism of groups
%(10)%%%%%%%%%%%%%%%%%%%%%%%%%%%%%%%%%%%%%%%%%%%%%%%%%%%%%%%%%%%%%%%%%%%%%%%
\begin{equation}\label{j}
\iota:~O(2n+1,q)~\rightarrow~Sp(2n,q)~(\begin{bmatrix}
  A & B & 0 \\
  C & D & 0 \\
  g & h & 1 \\
\end{bmatrix}~\mapsto~\begin{bmatrix}
                        A & B \\
                        C & D \\
                      \end{bmatrix}).
\end{equation}
%%%%%%%%%%%%%%%%%%%%%%%%%%%%%%%%%%%%%%%%%%%%%%%%%%%%%%%%%%%%%%%%%%%%%%%%%%%%%
In particular, for any $w\in O(2n+1,q)$,

%(11)%%%%%%%%%%%%%%%%%%%%%%%%%%%%%%%%%%%%%%%%%%%%%%%%%%%%%%%%%%%%%%%%%%%%%%%%
\begin{equation}\label{k}
Trw=Tr \iota (w)+1.
\end{equation}
%%%%%%%%%%%%%%%%%%%%%%%%%%%%%%%%%%%%%%%%%%%%%%%%%%%%%%%%%%%%%%%%%%%%%%%%%%%%%

Let $P=~P(2n+1,q)$ be the maximal parabolic subgroup of $O(2n+1,q)$ given by
\begin{align*}
P(2n+1,q)=\left\{\begin{bmatrix}
              A & 0 & 0 \\
              0 & ^t A^{-1} & 0 \\
              0 & 0 & 1 \\
            \end{bmatrix}\begin{bmatrix}
                           1_{n} & B & 0 \\
                           0 & 1_{n} & 0 \\
                           0 & h & 1 \\
                         \end{bmatrix} \Bigg\vert
\begin{array}{c}
  A\in GL(n,q) \\
  B+{^th} h~is ~alternating
\end{array}
\right\}.
\end{align*}\\\\
The Bruhat decomposition of $O(2n+1,q)$ with respect to $P=P(2n+1,q)$
is
\begin{align*}
O(2n+1,q)=~\coprod_{r~=~0}^{n}P\sigma_{r}P,
\end{align*}
where
\begin{align*}
\sigma_{r}~=~\begin{bmatrix}
               0 & 0 & 1_{r} & 0 & 0 \\
               0 & 1_{n-r} & 0 & 0 & 0 \\
               1_{r} & 0 & 0 & 0 & 0 \\
               0 & 0 & 0 & 1_{n-r} & 0 \\
               0 & 0 & 0 & 0 & 1 \\
             \end{bmatrix}~\in~O(2n+1,q).
\end{align*}\\\\

The symplectic group $Sp(2n,q)$ over the field $\mathbb{F}_{q}$ is defined as:
\begin{align*}
Sp(2n,q)~=~\{w\in GL(2n,q)|~^t wJw~=~J\},
\end{align*}
with
\begin{align*}
J=~\begin{bmatrix}
     0 & 1_{n} \\
     1_{n} & 0 \\
   \end{bmatrix}.
\end{align*}
Let $P^{'}=P^{'}(2n,q)$ be the maximal parabolic subgroup of
$Sp(2n,q)$ defined by:\\
\begin{align*}
P^{'}(2n,q)=~\left\{\begin{bmatrix}
                 A & 0 \\
                 0 & ^t A^{-1} \\
               \end{bmatrix}\begin{bmatrix}
                              1_{n} & B \\
                              0 & 1_{n} \\
                            \end{bmatrix}\Bigg|
                            ~A\in GL(n,q), ^t B=B\right\}.
\end{align*}\\
Then, with respect to $P^{'}=~P^{'}(2n,q)$, the Bruhat decomposition of
$Sp(2n,q)$ is given by
\begin{align*}
Sp(2n,q)=~\coprod_{r~=~0}^{n}P^{'}\sigma_{r}^{'}P^{'},
\end{align*}
where
\begin{align*}
\sigma_{r}^{'}=~\begin{bmatrix}
                  0 & 0 & 1_{r} & 0 \\
                  0 & 1_{n-r} & 0 & 0 \\
                  1_{r} & 0 & 0 & 0 \\
                  0 & 0 & 0 & 1_{n-r} \\
                \end{bmatrix}~\in~Sp(2n,q).
\end{align*}\\

Put, for each $r$ with $0 \leq r \leq n$,
\begin{gather*}
A_{r}=\{w\in P(2n+1,q)~|~\sigma_{r} w \sigma_{r}^{-1} \in P(2n+1,q)\},\\
A_{r}^{'}=\{w\in P^{'}(2n,q)~|~\sigma_{r}^{'}w(\sigma_{r}^{'})^{-1}\in P^{'}(2n,q)\}.
\end{gather*}
Expressing as the disjoint union of right cosets of maximal parabolic subgroups, the double cosets
$P\sigma_{r}P$ and $P^{'}\sigma_{r}^{'}P^{'}$ can be written respectively as

%(12)%%%%%%%%%%%%%%%%%%%%%%%%%%%%%%%%%%%%%%%%%%%%%%%%%%%%%%%%%%%%%%%%%%%%%%%%
\begin{equation}\label{l}
P\sigma_{r}P=~P\sigma_{r}(A_{r}\setminus P),
\end{equation}
%%%%%%%%%%%%%%%%%%%%%%%%%%%%%%%%%%%%%%%%%%%%%%%%%%%%%%%%%%%%%%%%%%%%%%%%%%%%%
%(13)%%%%%%%%%%%%%%%%%%%%%%%%%%%%%%%%%%%%%%%%%%%%%%%%%%%%%%%%%%%%%%%%%%%%%%%%
\begin{equation}\label{m}
P^{'}\sigma_{r}^{'}P^{'}=~P^{'}\sigma_{r}^{'}(A_{r}^{'}\setminus P^{'}).
\end{equation}\\
%%%%%%%%%%%%%%%%%%%%%%%%%%%%%%%%%%%%%%%%%%%%%%%%%%%%%%%%%%%%%%%%%%%%%%%%%%%%%

The order of the general linear group $GL(n,q)$ is given by

\begin{align*}
g_{n}=~\prod_{j=0}^{n-1}(q^{n}-q^{j})=~q^{\binom{n}{2}}\prod_{j=1}^{n}(q^{j}-1).
\end{align*}
For integers $n$,$r$ with $0\leq r\leq n$, the $q$-binomial coefficients are defined as:

\begin{align*}
\begin{bmatrix}
  n \\
  r \\
\end{bmatrix}_{q}=\prod_{j=0}^{r-1}(q^{n-j}-1)/(q^{r-j}-1).
\end{align*}

The following results follow either from \cite{D5} or from \cite{D1}
 plus the observation that under the isomorphism $\iota$ in (\ref{j})
 $P$, $A_{r}$, $\sigma_{r}$ are respectively mapped onto $P,^{'}A_{r}^{'},\sigma_{r}^{'}$:

%(14)%%%%%%%%%%%%%%%%%%%%%%%%%%%%%%%%%%%%%%%%%%%%%%%%%%%%%%%%%%%%%%%%%%%%%%%%%%%%%%%
\begin{equation}\label{n}
\begin{split}
&|A_{r}|=|A_{r}^{'}|=~g_{r}g_{n-r}q^{\binom{n+1}{2}}q^{r(2n-3r-1)/2},\\
&|P(2n+1,q)|=~|P^{'}(2n,q)|=q^{\binom{n+1}{2}}g_{n},\\
&|A_{r}\setminus P(2n+1,q)|=~|A_{r}^{'}\setminus P^{'}(2n,q)|=q^{\binom{r+1}{2}}\begin{bmatrix}
   n \\
   r \\
 \end{bmatrix}_{q},\\
&|P(2n+1,q)\sigma_{r}P(2n+1,q)|=~|P'(2n,q)\sigma_{r}^{'}P^{'}(2n,q)|\\
&\qquad \qquad \qquad\qquad\qquad\qquad =q^{n^{2}}\begin{bmatrix}
            n \\
            r \\
          \end{bmatrix}_{q}q^{\binom{r}{2}}q^{r}\prod_{j=1}^{n}(q^{j}-1)\\
&\qquad \qquad\qquad\qquad\qquad\qquad (=|P(2n+1,q)|^{2}|A_{r}|^{-1}\\
&\qquad \qquad\qquad\qquad\qquad\quad\quad =|P^{'}(2n+q)|^{2}|A_{r}^{'}|^{-1}).
\end{split}
\end{equation}
%%%%%%%%%%%%%%%%%%%%%%%%%%%%%%%%%%%%%%%%%%%%%%%%%%%%%%%%%%%%%%%%%%%%%%%%%%%%%%%%%%%%%%%%%%%%%%

In particular, with
\begin{align*}
DC(n,q)=~P(2n+1,q)\sigma_{n-1}P(2n+1,q),
\end{align*}

%(15)%%%%%%%%%%%%%%%%%%%%%%%%%%%%%%%%%%%%%%%%%%%%%%%%%%%%%%%%%%%%%%%%%%%%%%%%%%%%%%%%%%%%%%%%%%%
\begin{equation}\label{o}
|DC(n,q)|=~q^{\frac{1}{2}n(3n-1)}\begin{bmatrix}
                                   n \\
                                   1 \\
                                 \end{bmatrix}_{q}
\prod_{j=1}^{n}(q^{j}-1)=~A(n,q)B(n,q)~(cf.~ (\ref{d}),(\ref{e})).
\end{equation}\\\\
%%%%%%%%%%%%%%%%%%%%%%%%%%%%%%%%%%%%%%%%%%%%%%%%%%%%%%%%%%%%%%%%%%%%%%%%%%%%%%%%%%%%%%%%%%%%%%

%%%%%%%%%%%%%%%%%%%%%%%%%%%%%%%%%%%%%%%%%%%%%%%%%%%%%%%%%%%%%%%%%%%%%%%%
%%%%%%%%%%%%%%%%%%%%%%%%%%%%%%%%%%%%%%%%%%%%%%%%%%%%%%%%%%%%%%%%%%%%%%%%
\section{Exponential sums over double cosets of $O(2n+1,q)$}
%%%%%%%%%%%%%%%%%%%%%%%%%%%%%%%%%%%%%%%%%%%%%%%%%%%%%%%%%%%%%%%%%%%%%%%%
%%%%%%%%%%%%%%%%%%%%%%%%%%%%%%%%%%%%%%%%%%%%%%%%%%%%%%%%%%%%%%%%%%%%%%%%
The following notations will be employed throughout this paper.
\begin{align*}
\begin{split}
tr(x)&=x+x^2+\cdots+x^{2^{r-1}}~the~trace~function~\mathbb{F}_q\rightarrow\mathbb{F}_2,\\
 \lambda(x)&=(-1)^{tr(x)}~the~canonical~additive~character~of~\mathbb{F}_{q}.
\end{split}
\end{align*}
Then any nontrivial additive character $\psi$ of $\mathbb{F}_q$ is
given by $\psi(x)=\lambda(ax)$, for a unique
$a\in\mathbb{F}_{q}^{*}.$\\

For any nontrivial additive character $\psi$ of $\mathbb{F}_q$ and
$a\in\mathbb{F}_{q}^{*}$, the Kloosterman sum $K_{GL(t,q)}(\psi;a)$
for $GL(t,q)$ is defined as
\begin{align*}
K_{GL(t,q)}(\psi;a)=\sum_{w\in GL(t,q)}\psi(Trw+aTrw^{-1}).
\end{align*}
Notice that, for $t=1$, $K_{GL(1,q)}(\psi;a)$ denotes the Kloosterman sum
$K(\psi;a)$.\\
In \cite{D1}, it is shown that $K_{GL(t,q)}(\psi;a)$ satisfies the
following recursive relation: for integers $t\geq2$,
$a\in\mathbb{F}_{q}^{*}$,

\begin{align*}
K_{GL(t,q)}(\psi;a)=q^{t-1}K_{GL(t-1,q)}(\psi;a)K(\psi;a)+q^{2t-2}(q^{t-1}-1)K_{GL(t-2,q)}(\psi;a),
\end{align*}\\
where we understand that $K_{GL(0,q)}(\psi,a)=1$.\\

From \cite{D1} and \cite{D5}, we have (cf. (\ref{k})-(\ref{n})):

%{16}%%%%%%%%%%%%%%%%%%%%%%%%%%%%%%%%%%%%%%%%%%%%%%%%%%%%%%%%%%%%%%%%%%%%%%%%%%%%%%%%%%%
\begin{equation}\label{p}
\begin{split}
\sum_{w\in P\sigma_{r}P}&\psi(Trw)\\
&=|A_{r}\setminus P|\sum_{w\in P}\psi(Trw\sigma_{r})\\
&=\psi(1)|A_{r}^{'}\setminus P^{'}|\sum_{w\in P}\psi(Tr\iota(w)\sigma_{r}^{'})\\
&=\psi(1)|A_{r}^{'}\setminus P^{'}|\sum_{w\in P^{'}}\psi(Tr w\sigma_{r}^{'})\\
&(=\psi(1)\sum_{w\in P^{'}\sigma_{r}^{'}P^{'}}\psi(Trw))\\
&=\psi(1)q^{\binom{n+1}{2}}|A_{r}^{'}\setminus P^{'}|q^{r(n-r)}a_{r}K_{GL(n-r,q)}(\psi;1).
\end{split}
\end{equation}\\
%%%%%%%%%%%%%%%%%%%%%%%%%%%%%%%%%%%%%%%%%%%%%%%%%%%%%%%%%%%%%%%%%%%%%%%%%%%%%%%%%%%%%%

Here $\psi$ is any nontrivial additive character of $\mathbb{F}_{q}$, $a_0=1$, and,
for $r\in \mathbb{Z}_{>0}$, $a_{r}$ denotes the number of all $r\times r$ nonsingular
alternating matrices over $\mathbb{F}_{q}$, which is given by

%(17)%%%%%%%%%%%%%%%%%%%%%%%%%%%%%%%%%%%%%%%%%%%%%%%%%%%%%%%%%%%%%%%%%%%%%%%%%%%%%%%%%%%%%%%%%%%%
\begin{equation}\label{q}
a_{r}=\left\{
        \begin{array}{ll}
          0, & \hbox{if $r$ is odd,} \\
          q^{\frac{r}{2}(\frac{r}{2}-1)}\prod_{j=1}^{\frac{r}{2}}(q^{2j-1}-1), & \hbox{if $r$ is even}
        \end{array}
      \right.
\end{equation}
%%%%%%%%%%%%%%%%%%%%%%%%%%%%%%%%%%%%%%%%%%%%%%%%%%%%%%%%%%%%%%%%%%%%%%%%%%%%%%%%%%%%%%%%%%%%%%%
(cf.\cite{D1}, Proposition \ref{K}).\\
Thus we see from (\ref{n}), (\ref{p}), and (\ref{q}) that, for each $r$ with
$0\leq r\leq n$,

%(18)%%%%%%%%%%%%%%%%%%%%%%%%%%%%%%%%%%%%%%%%%%%%%%%%%%%%%%%%%%%%%%%%%%%%%%%%%%%%%%%%%%%%
\begin{equation}\label{r}
\sum_{w\in P\sigma_{r}P}\psi(Trw)=
\left\{
                                    \begin{array}{ll}
                                      0, & \hbox{if $r$ is odd,} \\
                                      \psi(1)q^{\binom{n+1}{2}}q^{rn-\frac{1}{4}r^2}\begin{bmatrix}
                        n \\
                        r \\
                      \end{bmatrix}_{q}\\
                                      \qquad\qquad\times\prod_{j=1}^{r/2}(q^{2j-1}-1)
K_{GL(n-r,q)}(\psi;1), & \hbox{if $r$ is even.} \\
                                    \end{array}
                                  \right.
\end{equation}\\
%%%%%%%%%%%%%%%%%%%%%%%%%%%%%%%%%%%%%%%%%%%%%%%%%%%%%%%%%%%%%%%%%%%%%%%%%%%%%%%%%%%%%%%

For our purposes, we need only one infinite family of exponential sums in (\ref{r}) over
$P(2n+1,q)\sigma_{n-1}P(2n+1,q)=~DC(n,q)$, for $n=1,3,5,\cdots$.
So we state them separately as a theorem.

%Theorem 2%%%%%%%%%%%%%%%%%%%%%%%%%%%%%%%%%%%%%%%%%%%%%%%%%%%%%%%%%%%%%%%%%%%%%%%%%
\begin{theorem}\label{B}
Let $\psi$ be any nontrivial additive character of $\mathbb{F}_q$.  Then in the
notation of $($\ref{d}$)$, we have
%(19)%%%%%%%%%%%%%%%%%%%%%%%%%%%%%%%%%%%%%%%%%%%%%%%%%%%%%%%%%%%%%%%%%%%%%%%%%%%%%%%%
\begin{equation}\label{s}
\sum_{w\in DC(n,q)}^{}\psi(Trw)=~\psi(1)A(n,q)K(\psi;1),for~ n=1,3,5,\cdots.
\end{equation}
%%%%%%%%%%%%%%%%%%%%%%%%%%%%%%%%%%%%%%%%%%%%%%%%%%%%%%%%%%%%%%%%%%%%%%%%%%%%%%%%%%%
\end{theorem}
%%%%%%%%%%%%%%%%%%%%%%%%%%%%%%%%%%%%%%%%%%%%%%%%%%%%%%%%%%%%%%%%%%%%%%%%%%%%%%%%%%%

%Proposition 3%%%%%%%%%%%%%%%%%%%%%%%%%%%%%%%%%%%%%%%%%%%%%%%%%%%%%%%%%%%%%%%%%%%%%
\begin{proposition}\label{C}$($\cite{D2}$)$
For $n=2^{s}(s\in \mathbb{Z}_{\geq 0})$, and $\lambda$ the nontrivial additive character
of $\mathbb{F}_q$,
\begin{align*}
K(\lambda;a^{n})=K(\lambda;a).
\end{align*}
\end{proposition}
%%%%%%%%%%%%%%%%%%%%%%%%%%%%%%%%%%%%%%%%%%%%%%%%%%%%%%%%%%%%%%%%%%%%%%%%%%%%%%%%%%%

The next corollary follows from Theorem $\ref{B}$, Proposition $\ref{C}$ and a simple
change of variables.

%Corollary 4%%%%%%%%%%%%%%%%%%%%%%%%%%%%%%%%%%%%%%%%%%%%%%%%%%%%%%%%%%%%%%%%%%%%%%%
\begin{corollary}\label{D}
Let $\lambda$ be the canonical additive character of $\mathbb{F}_q$, and let
$a\in \mathbb{F}_q^{*}$.  Then we have
%20%%%%%%%%%%%%%%%%%%%%%%%%%%%%%%%%%%%%%%%%%%%%%%%%%%%%%%%%%%%%%%%%%%%%%%%%%%%%%%%%
\begin{equation}\label{t}
\sum_{w\in DC(n,q)}\lambda(aTrw)=~\lambda(a)A(n,q)K(\lambda;a), for~
n=1,3,5,\cdots
\end{equation}
%%%%%%%%%%%%%%%%%%%%%%%%%%%%%%%%%%%%%%%%%%%%%%%%%%%%%%%%%%%%%%%%%%%%%%%%%%%%%%%%%%%
\end{corollary}
%%%%%%%%%%%%%%%%%%%%%%%%%%%%%%%%%%%%%%%%%%%%%%%%%%%%%%%%%%%%%%%%%%%%%%%%%%%%%%%%%%%
\noindent(cf. (\ref{d})).

%Proposition 5%%%%%%%%%%%%%%%%%%%%%%%%%%%%%%%%%%%%%%%%%%%%%%%%%%%%%%%%%%%%%%%%%%%%%
\begin{proposition}\label{E}$($\cite{D2}$)$
Let $\lambda$ be  the canonical additive character of $\mathbb{F}_q$,
$\beta\in\mathbb{F}_q$.  Then
%(21)%%%%%%%%%%%%%%%%%%%%%%%%%%%%%%%%%%%%%%%%%%%%%%%%%%%%%%%%%%%%%%%%%%%%%%%%%%%%%%%%
\begin{equation}\label{u}
\sum_{a\in \mathbb{F}_q^{*}}\lambda(-\alpha\beta)K(\lambda;a)
=\left\{
   \begin{array}{ll}
     q\lambda(\beta^{-1})+1, & \hbox{if $\beta$$\neq$0,} \\
     1, & \hbox{$if~\beta$=0.}
   \end{array}
 \right.
\end{equation}
%%%%%%%%%%%%%%%%%%%%%%%%%%%%%%%%%%%%%%%%%%%%%%%%%%%%%%%%%%%%%%%%%%%%%%%%%%%%%%%%%%%
\end{proposition}
%%%%%%%%%%%%%%%%%%%%%%%%%%%%%%%%%%%%%%%%%%%%%%%%%%%%%%%%%%%%%%%%%%%%%%%%%%%%%%%%%%%

For any integer $r$ with $0\leq r\leq n$, and each $\beta\in\mathbb{F}_q$, we let

\begin{align*}
N_{P\sigma_{r}P}(\beta)=|\{w\in P\sigma_{r}P|Trw=~\beta\}|.
\end{align*}\\
Then it is easy to see that

%(22)%%%%%%%%%%%%%%%%%%%%%%%%%%%%%%%%%%%%%%%%%%%%%%%%%%%%%%%%%%%%%%%%%%%%%%%%%%%%
\begin{equation}\label{v}
qN_{P\sigma_{r}P}(\beta)=|P\sigma_{r}P|+\sum_{a\in\mathbb{F}_{q}^{*}}\lambda(-\alpha\beta)
\sum_{w\in P\sigma_{r}P}\lambda(aTrw).
\end{equation}
%%%%%%%%%%%%%%%%%%%%%%%%%%%%%%%%%%%%%%%%%%%%%%%%%%%%%%%%%%%%%%%%%%%%%%%%%%%%%%%
For brevity, we write

%(23)%%%%%%%%%%%%%%%%%%%%%%%%%%%%%%%%%%%%%%%%%%%%%%%%%%%%%%%%%%%%%%%%%%%%%%%%%%%%
\begin{equation}\label{w}
n(\beta)=N_{DC(n,q)}(\beta).
\end{equation}\\
%%%%%%%%%%%%%%%%%%%%%%%%%%%%%%%%%%%%%%%%%%%%%%%%%%%%%%%%%%%%%%%%%%%%%%%%%%%%%%%
Now, from (\ref{o}) and (\ref{t})-(\ref{v}), we have the following result.

%Proposition 6%%%%%%%%%%%%%%%%%%%%%%%%%%%%%%%%%%%%%%%%%%%%%%%%%%%%%%%%%%%%%%%%%
\begin{proposition}\label{F}
With the notations in (\ref{d}), (\ref{e}), and (\ref{w}), for $n=1,3,5,\cdots,$
%24%%%%%%%%%%%%%%%%%%%%%%%%%%%%%%%%%%%%%%%%%%%%%%%%%%%%%%%%%%%%%%%%%%%%%%%%%%%%
\begin{equation}\label{x}
n(\beta)=q^{-1}A(n,q)B(n,q)+q^{-1}A(n,q)\times\left\{
                                                \begin{array}{ll}
                                                  1, & \hbox{$\beta=1$,} \\
                                                  q+1, & \hbox{$tr(\beta-1)^{-1}=0$,} \\
                                                  -q+1, & \hbox{$tr(\beta-1)^{-1}=1$.}
                                                \end{array}
                                              \right.
\end{equation}
%%%%%%%%%%%%%%%%%%%%%%%%%%%%%%%%%%%%%%%%%%%%%%%%%%%%%%%%%%%%%%%%%%%%%%%%%%%%%%%
\end{proposition}
%%%%%%%%%%%%%%%%%%%%%%%%%%%%%%%%%%%%%%%%%%%%%%%%%%%%%%%%%%%%%%%%%%%%%%%%%%%%%%%

%Corollary 7%%%%%%%%%%%%%%%%%%%%%%%%%%%%%%%%%%%%%%%%%%%%%%%%%%%%%%%%%%%%%%%%%%%
\begin{corollary}\label{G}
For each odd $n\geq3$, with all $q$, $n(\beta)>0$, for all $\beta$; for
$n=1$, with all $q$,
%25%%%%%%%%%%%%%%%%%%%%%%%%%%%%%%%%%%%%%%%%%%%%%%%%%%%%%%%%%%%%%%%%%%%%%%%%%%%%
\begin{equation}\label{y}
n(\beta)=\left\{
 \begin{array}{ll}
q, & \hbox{$\beta=1$,} \\
 2q, & \hbox{$tr(\beta-1)^{-1}=0$,} \\
 0, & \hbox{$tr(\beta-1)^{-1}=1$.}
 \end{array}
 \right.
\end{equation}
%%%%%%%%%%%%%%%%%%%%%%%%%%%%%%%%%%%%%%%%%%%%%%%%%%%%%%%%%%%%%%%%%%%%%%%%%%%%%%%
\end{corollary}
%%%%%%%%%%%%%%%%%%%%%%%%%%%%%%%%%%%%%%%%%%%%%%%%%%%%%%%%%%%%%%%%%%%%%%%%%%%%%%%

\begin{proof}
$n=1$ case follows directly from (\ref{x}).  Let $n\geq 3$ be odd.  Then, from
(\ref{x}), we see that, for any $\beta$,
\begin{align*}
\begin{split}
n(\beta)&\geq q^{\frac{1}{2}(3n^{2}-n-2)}(q^{n}-1)\prod_{j=2}^{n}(q^{j}-1)-q^{\frac{5}{4}(n^2-1)}\prod_{j=1}^{(n+1)/2}
(q^{2j-1}-1)\\
&>q^{\frac{1}{2}(3n^{2}-n-2)}(q^{n}-1)\prod_{j=2}^{n}(q^{j}-1)-q^{\frac{5}{4}(n^2-1)}\prod_{j=1}^{(n+1)/2}q^{2j-1}\\
&=q^{\frac{1}{2}(3n^{2}-n-2)}\{(q^{n}-1)(\prod_{j=2}^{n}(q^{j}-1)-1)-1\}>0.
\end{split}
\end{align*}
\end{proof}
%%%%%%%%%%%%%%%%%%%%%%%%%%%%%%%%%%%%%%%%%%%%%%%%%%%%%%%%%%%%%%%%%%%%%%%%%%%%%%%

%%%%%%%%%%%%%%%%%%%%%%%%%%%%%%%%%%%%%%%%%%%%%%%%%%%%%%%%%%%%%%%%%%%%%%%%
%%%%%%%%%%%%%%%%%%%%%%%%%%%%%%%%%%%%%%%%%%%%%%%%%%%%%%%%%%%%%%%%%%%%%%%%
\section{Construction of codes}
%%%%%%%%%%%%%%%%%%%%%%%%%%%%%%%%%%%%%%%%%%%%%%%%%%%%%%%%%%%%%%%%%%%%%%%%
%%%%%%%%%%%%%%%%%%%%%%%%%%%%%%%%%%%%%%%%%%%%%%%%%%%%%%%%%%%%%%%%%%%%%%%%
Let
%(26)%%%%%%%%%%%%%%%%%%%%%%%%%%%%%%%%%%%%%%%%%%%%%%%%%%%%%%%%%%%%%%%%%%%
\begin{equation}\label{z}
N(n,q)=|DC(n,q)|=A(n,q)B(n,q), ~for~n=1,3,5,\cdots
\end{equation}
%%%%%%%%%%%%%%%%%%%%%%%%%%%%%%%%%%%%%%%%%%%%%%%%%%%%%%%%%%%%%%%%%%%%%%%%
(cf. (\ref{d}), (\ref{e}), (\ref{o})).

Here we will construct one infinite family of binary linear codes
$C(DC(n,q))$ of length $N(n,q)$ for all positive odd integers $n$ and all $q$,
associated with the double cosets $DC(n,q)$.

Let $g_1,g_2,\cdots,g_{N(n,q)}$ be a fixed ordering of the elements in
$DC(n,q)$ $(n=1,3,5,\cdots).$  Then we put
\begin{align*}
v(n,q)=(Trg_1,Trg_2,\cdots,Trg_{N(n,q)})\in\mathbb{F}_{q}^{N(n,q)},~for~n=1,3,5,\cdots.
\end{align*}\\
Now, the binary linear code $C(DC(n,q))$ is defined as:

%27%%%%%%%%%%%%%%%%%%%%%%%%%%%%%%%%%%%%%%%%%%%%%%%%%%%%%%%%%%%%%%%%%%%%%%%%%%%%%
\begin{equation}\label{a1}
C(DC(n,q))=\{u\in\mathbb{F}_{2}^{N(n,q)}|u\cdot v(n,q)=0\},~for~n=1,3,5,\cdots,
\end{equation}
%%%%%%%%%%%%%%%%%%%%%%%%%%%%%%%%%%%%%%%%%%%%%%%%%%%%%%%%%%%%%%%%%%%%%%%%%%%%%%%%
where the dot denotes the usual inner product in $\mathbb{F}_{2}^{N(n,q)}.$

The following Delsarte's theorem is well-known.

%Theorem 8%%%%%%%%%%%%%%%%%%%%%%%%%%%%%%%%%%%%%%%%%%%%%%%%%%%%%%%%%%%%%%%%%%%%%%%
\begin{theorem}\label{H}$($\cite{L1}$)$
Let $B$ be a linear code over $\mathbb{F}_{q}.$  Then
\begin{align*}
(B|_{\mathbb{F}_{2}})^{\bot}=tr(B^{\bot}).
\end{align*}
\end{theorem}
%%%%%%%%%%%%%%%%%%%%%%%%%%%%%%%%%%%%%%%%%%%%%%%%%%%%%%%%%%%%%%%%%%%%%%%%%%%%%%%%%
\noindent
In view of this theorem, the dual $C(DC(n,q))^{\bot}$ of the code $C(DC(n,q))$
is given by

%28%%%%%%%%%%%%%%%%%%%%%%%%%%%%%%%%%%%%%%%%%%%%%%%%%%%%%%%%%%%%%%%%%%%%%%%%%%%%%%
\begin{equation}\label{b1}
\begin{split}
C(DC(&n,q))^{\bot}\\
&=~\{c(a)=c(a;n,q)=(tr(aTrg_1)),\cdots,tr(aTrg_{N(n,q)})|a\in\mathbb{F}_{q}\}
\end{split}
\end{equation}
%%%%%%%%%%%%%%%%%%%%%%%%%%%%%%%%%%%%%%%%%%%%%%%%%%%%%%%%%%%%%%%%%%%%%%%%%%%%%%%%%
($n=1,3,5,\cdots$).\\

Let $\mathbb{F}_{2}^{+}$, $\mathbb{F}_{q}^{+}$ denote the additive groups of the fields
$\mathbb{F}_{2}$, $\mathbb{F}_{q}$, respectively. Then we have the following exact
sequence of groups:

\begin{align*}
0\rightarrow\mathbb{F}_{2}^{+}\rightarrow\mathbb{F}_{q}^{+}\rightarrow
\Theta(\mathbb{F}_{q})\rightarrow0,
\end{align*}
where the first map is the inclusion and the second one is  the Artin-Schreier
operator in characteristic two given by $\Theta(x)=x^2+x.$  So
\begin{align*}
\Theta(\mathbb{F}_{q})=\{\alpha^2+\alpha|\alpha\in\mathbb{F}_{q}\}, ~and \quad[\mathbb{F}_{q}^{+} :\Theta(\mathbb{F}_{q})]=2.
\end{align*}

%Theorem 9%%%%%%%%%%%%%%%%%%%%%%%%%%%%%%%%%%%%%%%%%%%%%%%%%%%%%%%%%%%%%%%
\begin{theorem}\label{I}$($\cite{D2}$)$
Let $\lambda$ be the canonical additive character of $\mathbb{F}_{q}$,
and let $\beta\in\mathbb{F}_{q}^{*}$.  Then
%29%%%%%%%%%%%%%%%%%%%%%%%%%%%%%%%%%%%%%%%%%%%%%%%%%%%%%%%%%%%%%%%%%%%%%%%
\begin{equation}\label{c1}
\sum_{\alpha\in\mathbb{F}_{q}-\{0,1\}}^{}\lambda(\frac{\beta}{\alpha^2+\alpha})
=K(\lambda;\beta)-1.
\end{equation}
%%%%%%%%%%%%%%%%%%%%%%%%%%%%%%%%%%%%%%%%%%%%%%%%%%%%%%%%%%%%%%%%%%%%%%%%%%
\end{theorem}
%%%%%%%%%%%%%%%%%%%%%%%%%%%%%%%%%%%%%%%%%%%%%%%%%%%%%%%%%%%%%%%%%%%%%%%%%%

%Theorem 10%%%%%%%%%%%%%%%%%%%%%%%%%%%%%%%%%%%%%%%%%%%%%%%%%%%%%%%%%%%%%%%
\begin{theorem}\label{J}
The map $\mathbb{F}_{q}\rightarrow C(DC(n,q))^{\bot}(a\mapsto c(a))$ is an
$\mathbb{F}_{2}$-linear isomorphism for each odd integer $n\geq1$and all
$q$, except for $n=1$ and $q=4$.
\end{theorem}
%%%%%%%%%%%%%%%%%%%%%%%%%%%%%%%%%%%%%%%%%%%%%%%%%%%%%%%%%%%%%%%%%%%%%%%%%%

\begin{proof}
The map is clearly $\mathbb{F}_{2}$-linear and surjective. Let $a$ be in the
kernel of map. Then $tr(aTrg)=0$, for all $g\in DC(n,q)$.  If $n\geq3$
is odd, then, by Corollary \ref{G}, $Tr: DC(n,q)\rightarrow\mathbb{F}_{q}$
 is surjective and hence $tr(a\alpha)=0$, for all $\alpha\in\mathbb{F}_{q}$.
This implies that $a=0$, since otherwise $tr:\mathbb{F}_{q}\rightarrow\mathbb{F}_{2}$
would be the zero map. Now, assume that $n=1$.  Then, by (\ref{y}), $tr(a\beta)=0$,
for all $\beta\neq1$, with $tr((\beta-1)^{-1})=0$.
Hilbert's theorem 90 says that $tr(\gamma)=0\Leftrightarrow\gamma=\alpha^2+\alpha$,
for some $\alpha\in\mathbb{F}_{q}$.  This implies that
$\lambda(a)\sum_{\alpha\in\mathbb{F}_{q}-\{0,1\}}\lambda(\frac{a}{\alpha^2+\alpha})=
q-2$.  If $a\neq0,$ then, invoking (\ref{c1}) and the Weil bound (\ref{a}), we
would have
\begin{align*}
q-2=\lambda(a)\sum_{\alpha\in\mathbb{F}_{q}-\{0,1\}}\lambda(\frac{a}{\alpha^2+\alpha})=
\pm (K(\lambda;a)-1)\leq 2\sqrt{q}+1.
\end{align*}
For $q\geq 16$, this is impossible, since $x>2\sqrt{x}+3$, for $x\geq 16.$
On the other hand, for $q=2,4,8,$ one easily checks from (\ref{y}) that the kernel is trivial for
$q=2,8$ and is $\mathbb{F}_{2}$, for $q=4$.
\end{proof}
%%%%%%%%%%%%%%%%%%%%%%%%%%%%%%%%%%%%%%%%%%%%%%%%%%%%%%%%%%%%%%%%%%%%%%%%%%%%%%%%%%%%%%

%%%%%%%%%%%%%%%%%%%%%%%%%%%%%%%%%%%%%%%%%%%%%%%%%%%%%%%%%%%%%%%%%%%%%%%%
%%%%%%%%%%%%%%%%%%%%%%%%%%%%%%%%%%%%%%%%%%%%%%%%%%%%%%%%%%%%%%%%%%%%%%%%
\section{Power moments of Kloosterman sums with trace one arguments}
%%%%%%%%%%%%%%%%%%%%%%%%%%%%%%%%%%%%%%%%%%%%%%%%%%%%%%%%%%%%%%%%%%%%%%%%
%%%%%%%%%%%%%%%%%%%%%%%%%%%%%%%%%%%%%%%%%%%%%%%%%%%%%%%%%%%%%%%%%%%%%%%%

Here we will be able to find, via Pless power moment identity, an infinite
family of recursive formulas generating the odd power moments of  Kloosterman
sums with trace one arguments over all $\mathbb{F}_{q}$
in terms of the frequencies of weights in
$C(DC(n,q))$ and $C(\widehat{DC}(n,q))$, respectively.

%Theorem 11%%%%%%%%%%%%%%%%%%%%%%%%%%%%%%%%%%%%%%%%%%%%%%%%%%%%%%%%%%%%%%%%%%
\begin{theorem}\label{K}$($Pless power moment identity, \cite{L1}$)$
 Let $B$ be an $q$-ary $[n,k]$ code, and let $B_i$(resp. $B_i^\bot$)
 denote the number of codewords of weight $i$ in $B$(resp. in
 $B^\bot$). Then, for $h=0,1,2,\cdots,$
%(30)%%%%%%%%%%%%%%%%%%%%%%%%%%%%%%%%%%%%%%%%%%%%%%%%%%%%%%%%%%%%%%%%%%%%%%%%
\begin{equation}\label{d1}
\sum_{j=0}^{n}j^hB_j=\sum_{j=0}^{min\{n,h\}}(-1)^jB_j^\bot
\sum_{t=j}^{h}t!S(h,t)q^{k-t}(q-1)^{t-j}{\binom{n-j}{n-t}},
\end{equation}
%%%%%%%%%%%%%%%%%%%%%%%%%%%%%%%%%%%%%%%%%%%%%%%%%%%%%%%%%%%%%%%%%%%%%%%%%%%%%
where $S(h,t)$ is the Stirling number of the second kind defined in
$($\ref{i}$)$.
\end{theorem}
%%%%%%%%%%%%%%%%%%%%%%%%%%%%%%%%%%%%%%%%%%%%%%%%%%%%%%%%%%%%%%%%%%%%%%%%%%%%%%

%Lemma 12%%%%%%%%%%%%%%%%%%%%%%%%%%%%%%%%%%%%%%%%%%%%%%%%%%%%%%%%%%%%%%%%%%%%%
\begin{lemma}\label{L}
Let $c(a)=(tr(aTrg_1),\cdots,tr(aTrg_{N(n,q)}))\in C(DC(n,q))^\bot$$($$n=1,3,5,\cdots$$)$,
for $a\in\mathbb{F}_{q}^{*}$.  Then the Hamming weight
$w(c(a))$ is expressed as follows:
%31%%%%%%%%%%%%%%%%%%%%%%%%%%%%%%%%%%%%%%%%%%%%%%%%%%%%%%%%%%%%%%%%%%%%%%%%%%%
\begin{equation}\label{e1}
w(c(a))=\frac{1}{2}A(n,q)(B(n,q)-\lambda(a)K(\lambda;a))(cf.~(\ref{d}),~(\ref{e})).
\end{equation}
%%%%%%%%%%%%%%%%%%%%%%%%%%%%%%%%%%%%%%%%%%%%%%%%%%%%%%%%%%%%%%%%%%%%%%%%%%%%%%
\end{lemma}

\begin{proof}
\begin{align*}
w(c(a))=\frac{1}{2}\sum_{j=1}^{N(n,q)}(1-(-1)^{tr(aTrg_{j})})
=\frac{1}{2}(N(n,q)-\sum_{w\in DC(n,q)}\lambda(aTrw)).
\end{align*}
Our result now follows from (\ref{t}) and (\ref{z}).
\end{proof}
%%%%%%%%%%%%%%%%%%%%%%%%%%%%%%%%%%%%%%%%%%%%%%%%%%%%%%%%%%%%%%%%%%%%%%%%%%%%%%

Let $u=(u_1,\cdots,u_{N_{N(n,q)}})\in\mathbb{F}_{2}^{N(n,q)}$, with $\nu_{\beta}$
1's in the coordinate palces where $Tr(g_j)=\beta$, for each $\beta\in\mathbb{F}_{q}$.
Then from the definition of the codes $C(DC(n,q))$ (cf.(\ref{a1}))
that $u$ is a codeword with weight $j$ if and only if
$\sum_{\beta\in\mathbb{F}_{q}}\nu_{\beta}=j$ and $\sum_{\beta\in\mathbb{F}_{q}}\nu
_{\beta}\beta=0$(an identity in $\mathbb{F}_{q}$).  As there are
$\prod_{\beta\in\mathbb{F}_{q}}\begin{pmatrix}
                                 n(\beta) \\
                                 \nu_{\beta} \\
                               \end{pmatrix}$
(cf. (\ref{w})) many such codewords with weight $j$, we obtain the following result.

%Proposition 13%%%%%%%%%%%%%%%%%%%%%%%%%%%%%%%%%%%%%%%%%%%%%%%%%%%%%%%%%%%%%%%%
\begin{proposition}\label{M}
Let $\{C_j(n,q)\}_{j=0}^{N(n,q)}$ be the weight distribution of
$C(DC(n,q))$ $(n=1,3,5,\cdots)$. Then
%32%%%%%%%%%%%%%%%%%%%%%%%%%%%%%%%%%%%%%%%%%%%%%%%%%%%%%%%%%%%%%%%%%%%%%%%%%%%%
\begin{equation}\label{f1}
C_j(n,q)=~\sum\prod_{\beta\in\mathbb{F}_{q}}\begin{pmatrix}
                                              n(\beta) \\
                                              \nu_{\beta} \\
                                            \end{pmatrix},
\end{equation}
%%%%%%%%%%%%%%%%%%%%%%%%%%%%%%%%%%%%%%%%%%%%%%%%%%%%%%%%%%%%%%%%%%%%%%%%%%%%%%%
where the sum is over all the sets of integers $\{\nu_{\beta}\}_{\beta\in\mathbb{F}_{q}}$
$(0\leq\nu_{\beta}\leq n(\beta))$, satisfying
%33%%%%%%%%%%%%%%%%%%%%%%%%%%%%%%%%%%%%%%%%%%%%%%%%%%%%%%%%%%%%%%%%%%%%%%%%%%%%
\begin{equation}\label{g1}
\sum_{\beta\in\mathbb{F}_{q}}\nu_{\beta}=j,~ and
\sum_{\beta\in\mathbb{F}_{q}}\nu_{\beta}\beta=0.
\end{equation}
%%%%%%%%%%%%%%%%%%%%%%%%%%%%%%%%%%%%%%%%%%%%%%%%%%%%%%%%%%%%%%%%%%%%%%%%%%%%%%%
\end{proposition}
%%%%%%%%%%%%%%%%%%%%%%%%%%%%%%%%%%%%%%%%%%%%%%%%%%%%%%%%%%%%%%%%%%%%%%%%%%%%%%%

%Corollary 14%%%%%%%%%%%%%%%%%%%%%%%%%%%%%%%%%%%%%%%%%%%%%%%%%%%%%%%%%%%%%%%%%%
\begin{corollary}\label{N}
Let $\{C_{j}(n,q)\}_{j=0}^{N(n,q)}$ $(n=1,3,5,\cdots)$ be as above.  Then
we have
\begin{align*}
C_j(n,q)=C_{N(n,q)-j}(n,q),~ for ~all~ j,~ with~ 0\leq j\leq N(n,q).
\end{align*}
\end{corollary}
%%%%%%%%%%%%%%%%%%%%%%%%%%%%%%%%%%%%%%%%%%%%%%%%%%%%%%%%%%%%%%%%%%%%%%%%%%%%%%%

\begin{proof}
Under the replacements $\nu_{\beta}\rightarrow n(\beta)-\nu_{\beta}$, for
each $\beta\in\mathbb{F}_{q}$, the first equation in $($\ref{g1}$)$ is changed
to $N(n,q)-j$, while the second one in there and the summand in
$($\ref{f1}$)$ is left unchanged. Here the second sum in $($\ref{g1}$)$
is left unchanged, since $\sum_{\beta\in\mathbb{F}_{q}}n(\beta)\beta=0$,
as one can see by using the explicit expressions of $n(\beta)$
in $($\ref{x}$)$ and $($\ref{y}$)$.
\end{proof}
%%%%%%%%%%%%%%%%%%%%%%%%%%%%%%%%%%%%%%%%%%%%%%%%%%%%%%%%%%%%%%%%%%%%%%%%%%%%%%%

The formula appearing in the next theorem and stated in (\ref{g})
follows from the formula in (\ref{f1}), using the explicit value of
$n(\beta)$ in (\ref{x}).

%Theorem 15%%%%%%%%%%%%%%%%%%%%%%%%%%%%%%%%%%%%%%%%%%%%%%%%%%%%%%%%%%%%%%%%%%%%
\begin{theorem}\label{O}
Let $\{C_{j}(n,q)\}_{j=0}^{N(n,q)}$ be the weight distribution of
$C(DC(n,q))$ $(n=1,3,5,\cdots)$. Then, for $j=0,\cdots,N(n,q)$,
\begin{align*}
\begin{split}
C_{j}(n,q)&=\sum\begin{pmatrix}
               q^{-1}A(n,q)(B(n,q)+1) \\
               \nu_1 \\
             \end{pmatrix}\\
&\times\prod_{tr(\beta-1)^{-1}=0}\begin{pmatrix}
                                   q^{-1}A(n,q)(B(n,q)+q+1) \\
                                   \nu_{\beta} \\
                                 \end{pmatrix}\\
&\qquad\qquad\qquad\times\prod_{tr(\beta-1)^{-1}=1}\begin{pmatrix}
                                   q^{-1}A(n,q)(B(n,q)-q+1) \\
                                   \nu_{\beta} \\
                                 \end{pmatrix},
\end{split}
\end{align*}
where the sum is over all the sets of nonnegative integers
$\{\nu_{\beta}\}_{\beta\in\mathbb{F}_{q}}$ satisfying
$\sum_{\beta\in\mathbb{F}_{q}}\nu_{\beta}=j$ and
$\sum_{\beta\in\mathbb{F}_{q}}\nu_{\beta}\beta=0.$
\end{theorem}
%%%%%%%%%%%%%%%%%%%%%%%%%%%%%%%%%%%%%%%%%%%%%%%%%%%%%%%%%%%%%%%%%%%%%%%%%%%%%%%

The recursive formula in the next theorem follows from the study of codes
associated with the double cosets $\widehat{DC}(n,q)=P^{'}(2n,q)\sigma
_{n-1}^{'}P^{'}(2n,q)$ of the symplectic group $Sp(2n,q)$.
It is slightly modified from its original version, which makes it more usable in below.

%Theorem 16%%%%%%%%%%%%%%%%%%%%%%%%%%%%%%%%%%%%%%%%%%%%%%%%%%%%%%%%%%%%%%%%%%%%%%
\begin{theorem}\label{P}$($\cite{D3}$)$
For each odd integer $n\geq 3$, with all $q$, or $n=1$, with $q\geq8$,
%34%%%%%%%%%%%%%%%%%%%%%%%%%%%%%%%%%%%%%%%%%%%%%%%%%%%%%%%%%%%%%%%%%%%%%%%%%%%%%%
\begin{equation}\label{h1}
\begin{split}
&\frac{1}{2^h}A(n,q)^h\sum_{l=0}^{h}(-1)^l\begin{pmatrix}
                                           h \\
                                           l \\
                                         \end{pmatrix}B(n,q)^{h-l}MK^l\\
&=q\sum_{j=0}^{min\{N(n,q),h\}}(-1)^j\widehat{C}_j(n,q)\sum_{t=j}^{h}t!
S(h,t)2^{-t}\begin{pmatrix}
              N(n,q)-j \\
              N(n,q)-t \\
            \end{pmatrix}(h=1,2,\cdots),
\end{split}
\end{equation}
%%%%%%%%%%%%%%%%%%%%%%%%%%%%%%%%%%%%%%%%%%%%%%%%%%%%%%%%%%%%%%%%%%%%%%%%%%%%%%%%%
where $N(n,q)=A(n,q)B(n,q)$, and $\{\widehat{C}_j(n,q)\}_{j=0}^{N(n,q)}$
is the weight distribution of $C(\widehat{DC}(n,q))$ given by
\begin{align*}
\begin{split}
\widehat{C}_{j}(n,q)&=\sum\begin{pmatrix}
                           q^{-1}A(n,q)(B(n,q)+1) \\
                           \nu_0 \\
                         \end{pmatrix}\\
&\times\prod_{tr(\beta^{-1})=0}\begin{pmatrix}
                                 q^{-1}A(n,q)(B(n,q)+q+1) \\
                                 \nu_{\beta} \\
                               \end{pmatrix}\\
&\qquad\qquad\qquad\times\prod_{tr(\beta^{-1})=1}\begin{pmatrix}
                                                   q^{-1}A(n,q)(B(n,q)-q+1) \\
                                                   \nu_{\beta} \\
                                                 \end{pmatrix}.
\end{split}
\end{align*}
\end{theorem}
%%%%%%%%%%%%%%%%%%%%%%%%%%%%%%%%%%%%%%%%%%%%%%%%%%%%%%%%%%%%%%%%%%%%%%%%%%%%%%%%%

Here the sum is over all the sets of nonnegative integers $\{\nu_{\beta}\}_{\beta\in\mathbb{F}_{q}}$
satisfying $\sum_{\beta\in\mathbb{F}_{q}}\nu_{\beta}=j$ and
$\sum_{\beta\in\mathbb{F}_{q}}\nu_{\beta}\beta=0$. In addition,
$S(h,t)$ is the Stirling number of the second kind as in
(\ref{i}).\\

From now on, we will assume that $n$ is any odd integer $\geq3$, with all $q$, or
$n=1$, with $q\geq8$.  Under these assumptions, each codeword in
$C(DC(n,q))^{\bot}$ can be written as $c(a)$, for a unique $a\in\mathbb{F}_{q}$
(cf. Theorem \ref{J}, (\ref{b1})) and Theorem \ref{P} in the above can be applied.

  Now, we apply the Pless power moment identity in (\ref{d1}) to
$C(DC(n,q))^{\bot}$, in order to get the result in Theorem \ref{A} (cf.
(\ref{f})-(\ref{h})) about recursive formulas. In below, ``the sum over
$tra=0$(resp. $tra=1$)" will mean ``the sum over all nonzero $a\in\mathbb{F}_{q}^{*}$,
with $tr a=0$ (resp. $tr a=1$)."
The left hand side of that identity in (\ref{d1}) is equal to

\begin{align*}
\sum_{a\in\mathbb{F}_{q}^{*}}w(c(a))^h,
\end{align*}
with $w(c(a))$ given by (\ref{e1}).  We have

%35%%%%%%%%%%%%%%%%%%%%%%%%%%%%%%%%%%%%%%%%%%%%%%%%%%%%%%%%%%%%%%%%%%%%%%%%%%%%%
\begin{equation}\label{i1}
\begin{split}
&\sum_{a\in\mathbb{F}_{q}^{*}}w(c(a))^h=\frac{1}{2^h}A(n,q)^h
\sum_{a\in\mathbb{F}_{q}^{*}}(B(n,q)-\lambda(a)K(\lambda;a))^h\\
&=\frac{1}{2^h}A(n,q)^h\sum_{tra=0}(B(n,q)-K(\lambda;a))^h
+\frac{1}{2^h}A(n,q)^h\sum_{tra=1}(B(n,q)+K(\lambda;a))^h\\
&=\frac{1}{2^h}A(n,q)^h\sum_{tra=0}\sum_{l=0}^{h}(-1)^l\begin{pmatrix}
                                                        h \\
                                                        l \\
                                                      \end{pmatrix}
B(n,q)^{h-l}K(\lambda;a)^l\\
&\qquad\qquad +\frac{1}{2^h}A(n,q)^h\sum_{tra=1}\sum_{l=0}^{h}\begin{pmatrix}
                                                        h \\
                                                        l \\
                                                      \end{pmatrix}
B(n,q)^{h-l}K(\lambda;a)^l\\
&=\frac{1}{2^h}A(n,q)^h\sum_{l=0}^{h}(-1)^l\begin{pmatrix}
                                                        h \\
                                                        l \\
                                                      \end{pmatrix}
B(n,q)^{h-l}(MK^l-T_1K^l)(\psi=\lambda~case~of~(\ref{b}),~(\ref{c}))\\
&\qquad\qquad +\frac{1}{2^h}A(n,q)^h\sum_{l=0}^{h}\begin{pmatrix}
                                                        h \\
                                                        l \\
                                                      \end{pmatrix}
B(n,q)^{h-l}T_1K^l\\
&=\frac{1}{2^h}A(n,q)^h\sum_{l=0}^{h}(-1)^l\begin{pmatrix}
                                                        h \\
                                                        l \\
                                                      \end{pmatrix}
B(n,q)^{h-l}MK^l\\
&\qquad\qquad +2\frac{1}{2^h}A(n,q)^h
\sum_{0\leq l\leq h,~ l~ odd}\begin{pmatrix}
                                                        h \\
                                                        l \\
                                                      \end{pmatrix}
B(n,q)^{h-l}T_1K^l\\
&=q\sum_{j=0}^{min\{N(n,q),h\}}(-1)^j\widehat{C}_j(n,q)
\sum_{t=j}^{h}t!S(h,t)2^{-t}\begin{pmatrix}
                                                      N(n,q)-j \\
                                                      N(n,q)-t \\
                                                      \end{pmatrix}(cf. ~(\ref{h1}))\\
&\qquad\qquad +2\frac{1}{2^h}A(n,q)^h
\sum_{0\leq l\leq h,~ l~ odd}\begin{pmatrix}
                                                        h \\
                                                        l \\
                                                      \end{pmatrix}
B(n,q)^{h-l}T_1K^{l}.\\
\end{split}
\end{equation}
%%%%%%%%%%%%%%%%%%%%%%%%%%%%%%%%%%%%%%%%%%%%%%%%%%%%%%%%%%%%%%%%%%%%%%%%%%%%%%%%

On the other hand, the right hand side of the identity in (\ref{d1}) is
given by:

%36%%%%%%%%%%%%%%%%%%%%%%%%%%%%%%%%%%%%%%%%%%%%%%%%%%%%%%%%%%%%%%%%%%%%%%%%%%%%%
\begin{equation}\label{j1}
q\sum_{j=0}^{min\{N(n,q),h\}}(-1)^{j}C_j(n,q)
\sum_{t=j}^{h}t!S(h,t)2^{-t}\begin{pmatrix}
                                                      N(n,q)-j \\
                                                      N(n,q)-t \\
                                                      \end{pmatrix}.
\end{equation}
%%%%%%%%%%%%%%%%%%%%%%%%%%%%%%%%%%%%%%%%%%%%%%%%%%%%%%%%%%%%%%%%%%%%%%%%%%%%%%%%

In (\ref{j1}), one has to note that $dim_{\mathbb{F}_{2}}C(DC(n,q))^{\bot}=r$.
Our main result in (\ref{f}) now follows by equating (\ref{i1}) and (\ref{j1}).\\\\

\bibliographystyle{amsplain}

\end{document}